\documentclass[12pt]{amsart}
\usepackage{amssymb}
\usepackage{verbatim}
\usepackage{xspace}
\usepackage{vmargin}
\setpapersize{USletter}
\setmarginsrb{1.1in}{1.2in}{1.1in}{1.2in}{12 pt}{25 pt}{12 pt}{30 pt}

\theoremstyle{plain}
\newtheorem{theorem}{Theorem}

\theoremstyle{definition}

\theoremstyle{remark}
\newtheorem*{remark}{Remark}

\newcommand{\ds}{\displaystyle}

\newcommand{\cstar}{\ensuremath{\text{C}^{*}}-}
\renewcommand{\star}{\ensuremath{{}^{*}}\nobreakdash-\hspace{0 pt}}
\newcommand{\ol}[1]{\overline{#1}}
\newcommand{\sads}[1]{S_{#1}^* S_{#1}^{\vphantom{*}}}
\newcommand{\ssad}[1]{S_{#1}^{\vphantom{*}} S_{#1}^*}

\newcommand{\brk}{{\beta_{(r,k)}}} 
\newcommand{\Urk}{{U^{(r,k)}}}
\newcommand{\Irk}{{I_{(r,k)}}}
\newcommand{\Trk}{T_{(r,k)}}

\DeclareMathOperator{\dom}{dom}
\DeclareMathOperator{\ran}{ran}

\newcommand{\bbF}{\mathbb F}

\newcommand{\bbQ}{\mathbb Q}

\newcommand{\bbZ}{\mathbb Z}

     \newcommand{\sB}{\mathcal B}

     \newcommand{\sF}{\mathcal F}
     \newcommand{\sG}{\mathcal G}
     \newcommand{\sH}{\mathcal H}

     \newcommand{\sV}{\mathcal V}

\begin{document}
% topmatter
\title[Partial Crossed Product Presentations for
               $O_n$ and $M_k(O_n)$]{Partial 
    Crossed Product Presentations for\\
               $O_n$ and $M_k(O_n)$ Using Amenable Groups}
\author{Alan Hopenwasser}
\address{Department of Mathematics\\
        University of Alabama\\
        Tuscaloosa, AL 35487}
\email{ahopenwa@bama.ua.edu}

 \date{February 17, 2006}

% abstract
 \begin{abstract}
The Cuntz algebra $O_n$ is presented as a partial crossed product in which an amenable group partially acts on an abelian {\ensuremath{\text{C}^{*}}-}algebra.  The partial action is related to the Cuntz groupoid for $O_n$ and connections are made with non-self-adjoint subalgebras of $O_n$, particularly the Volterra nest subalgebra.  These ideas are also extended to the $M_k(O_n)$ context.
 \end{abstract}
\maketitle

\section{Introduction} \label{s:intro}

The purpose of this note is to present the Cuntz algebra,
$O_n$, as a partial crossed product by an amenable group acting on an abelian \cstar algebra; to connect this presentation with the standard presentation of $O_n$ as a groupoid \cstar algebra based on the Cuntz groupoid; and to indicate how the partial crossed product presentation ties in with certain non-self-adjoint subalgebras of $O_n$, most notably, the Volterra nest subalgebra.  In addition, we indicate how each $M_k (O_n)$ can also be written as a partial crossed product by the partial action of an amenable group on an abelian \cstar algebra.  We also connect the partial crossed product presentation with a groupoid isomorphic to the usual groupoid for $M_k(O_n)$ obtained by viewing $M_k(O_n)$ as a graph \cstar algebra.

There is nothing new about writing $O_n$ (or $M_k (O_n)$) as a crossed product or partial crossed product.  In \cite{MR1113572} and
\cite{MR1245830}, Spielberg shows that $O_n$ (as well as other graph \cstar algebras) can be presented as a crossed product in which the group is a free product of cyclic groups.  Quigg and Raeburn in
\cite{MR1452280} show that $O_n$ is a partial crossed product using the free group, $\bbF_n$, on $n$ generators.  Exel, Laca, and Quigg extend this to graph \cstar algebras (with finite graph), again using 
$\bbF_n$.  In all these presentations, the group used is not amenable.
The group appearing in our presentation of $O_n$ is a semi-direct product of the integers and the group $\bbQ_n$ of $n$-adic rationals, and so is amenable.  This immediately yields yet another proof that $O_n$ is nuclear.  The group used in the presentation of $M_k(O_n)$ is the direct product of the group $\bbQ_n \times_{\delta} \bbZ$
 and a cyclic group, and so is also amenable.
  There is nothing canonical about this presentation (or about the related groupoids); but they do display much about the internal structure of $O_n$ (and
$M_k (O_n)$) and as such should prove useful in the study of non-self-adjoint subalgebras of these \cstar algebras.

The author wishes to thank a number of people for helpful conversations on the subject matter of this note: these include Allan Donsig, Ruy Exel, Justin Peters, David Pitts, and Steve Power.  The author also thanks the referee for helpful comments which resulted in a more intuitive proof of Theorem
\ref{t:Onpca}.  The original proof was long, tedious, and not particularly illuminating.

\section{$O_n$ as a partial crossed product}  \label{s:On}

Throughout this section the positive integer $n \geq 2$ will be
fixed.  Let $X$  denote the Cantor set based on $[0,1]$ with
each $n$-adic rational $r$ replaced by a pair $r^-$, $r^+$.  (0 and 1
are excepted, of course, and $r^-$ is the immediate predecessor of
$r^+$.)  The group, $G$, which will partially act on the abelian \cstar
algebra $C(X)$ is the semi-direct product of the additive group 
$\bbQ_n$ of $n$-adic rationals and the integers.  The action,
$\delta$, of $\bbZ$
on $\bbQ_n$ is given by
\begin{displaymath}
  \delta_k(r) = \frac r{n^k}, \text{ for all } k \in \bbZ, r \in \bbQ_n.
\end{displaymath}
Multiplication and inverses in $G = \bbQ_n \times_{\delta} \bbZ$ are
given by
\begin{align*}
&(s,j)(r,k) = (\delta_j(r)+s,j+k)
=(\frac r{n^j}+s, j+k), \\
&(r,k)^{-1} = (-\delta_{-k}(r), -k)
=(-n^kr,-k). 
\end{align*}

Define a partial action $\beta$ of $G$ on $X$ by
\[
\beta_{(r,k)}(x) =\frac x{n^k}+r.
\]
The domain of $\beta_{(r,k)}$ is the set of all $x \in X$ for which 
$\frac x{n^k}+r \in X$. 
(If $x$ is $n$-adic, then
$\beta_{(r,k)}(x^-) = (\frac x{n^k}+r)^-$ and
$\beta_{(r,k)}(x^+) = (\frac x{n^k}+r)^+$.  However, normally
we drop the $+$ or $-$ superscript.)

As usual, the partial action on $X$ induces a partial action $\alpha$
on $C(X)$.  This is given by
\begin{displaymath}
  \alpha_g(f) = f \circ \beta_g^{-1}, \quad 
  f \in C_0(\dom \beta_g), \quad g \in G.
\end{displaymath}
While $f \circ \beta_g^{-1}$ is defined only on
$\ran \beta_g$, it is sometimes
helpful to interpret it to be defined on all of
$X$ by declaring that it has value 0 at points not in
$\ran \beta_g$.

The partial crossed product $C(X) \times_{\alpha} G$ is generated by
polynomials in the symbols $U^g$ with the coefficient of
$U^g$ always in 
$C_0(\ran \beta_g) = \ran \alpha_g$.  
We should point out that the symbols $U^g$ will not be elements
of $C(X) \times_{\alpha} G$; coefficients in
$C_0(\ran \beta_g) = \ran \alpha_g$ are mandatory.  However,
 these ideals have units, so the elements
$\chi_{\ran \beta_g} U^g$ do serve as substitutes for the
$U^g$.  Be warned, however, that 
$\chi_{\ran \beta_g} U^g \chi_{\ran \beta_h} U^h$ is not 
generally equal to $\chi_{\ran \beta_{gh}} U^{gh}$.  This is a consequence 
of the fact that $\beta_g \circ \beta_h$ is a (usually proper)
restriction of $\beta_{gh}$.
The primary
algebraic formulas which are needed to verify various claims in
the proof of Theorem \ref{t:Onpca}
 are the formulas for inversion and
multiplication on monomials:
  \begin{align}
\left(f U^g\right)^* 
&= \ol{\alpha_{g^{-1}}(f)}U^{g^{-1}} \label{inv}\\
&= \ol{f \circ \beta_{g}} U^{g^{-1}}\text{ and} \notag\\
eU^{g} fU^{h} 
&= \alpha_{g} (\alpha_{g^{-1}}(e)f)
  U^{gh} \label{mult}\\
&= e \cdot f \circ \beta_{g^{-1}} U^{gh},\notag
\end{align}
We shall use primarily  the variant with $\beta$ in the sequel.

\begin{remark}
The group $\bbQ_n$ embeds naturally in $G$.  If we restrict the partial action $\alpha$ to $\bbQ_n$, then $C(X) \times_{\alpha} \bbQ_n$ is an $n^\infty$-UHF algebra.  Indeed, viewed as a subalgebra of
$C(X) \times_{\alpha} G$ this will turn out to be (after the identification of $C(X) \times_{\alpha} G$ as $O_n$) the core UHF subalgebra of $O_n$.
This subalgebra appears in section \ref{subalg} where it plays a role in identifying some non-self-adjoint subalgebras of $O_n$ in the partial crossed product formulation.
\end{remark}

\begin{theorem} \label{t:Onpca}
 $C(X) \times_{\alpha} G$ is isomorphic to the Cuntz algebra $O_n$.
\end{theorem}

The proof of Theorem \ref{t:Onpca} will make use of the standard representation of $O_n$ as a concrete operator algebra acting on
$L^2[0,1]$.  In the discussion of this representation, all subintervals of $[0,1]$ will have $n$-adic endpoints; so take the word ``interval'' to mean ``interval with $n$-adic endpoints''.  If $\phi$ is an order preserving affine transformation (hereafter abbreviated ``opat'') from a subinterval $D$ of $[0,1]$ onto a subinerval $R$, then $\phi$ induces a partial isometry on $L^2[0,1]$ in the usual way ($f \mapsto f \circ \phi^{-1}$) with initial and final spaces identifiable with $L^2(D)$ and $L^2(R)$.  Denote this partial isometry by $T(\phi)$.  In particular, for each
$(r,k) \in G$ the formula for $\beta_{(r,k)}$ on the Cantor space $X$ given above (viz. $x \mapsto \frac x{n^k} +r$) gives an opat from a subinterval of $[0,1]$ to another subinterval.  Using this formula, we obtain a partial action of $G$ on $[0,1]$, which we shall also denote by $\beta$.  (The context will determine clearly which partial action is intended.)  The domain of $\beta_{(r,k)}$ is, of course, the largest interval in $[0,1]$ which is mapped again into $[0,1]$.

For each $i=1, \dots, n$, $\beta_{\left(\frac {i-1}n, 1\right)}$ is the opat from $[0,1]$ onto $[\frac {i-1}n, \frac in]$; the isometries
$T_i = T(\beta_{\left(\frac {i-1}n, 1\right)})$ have mutually orthogonal range projections which add to $I$.  The representation of $O_n$ which we use is the \cstar algebra generated by $T_1, \dots, T_n$.  For the rest of this section, $O_n$ refers to this algebra.  

Let $\sF$ denote the family of all opats $\phi$ which are restrictions of some $\beta_{(r,k)}$.  If $\phi$ and $\psi$ are two opats in $\sF$, then
$\phi \circ \psi$ denotes the opat whose domain is
$\{ x \in \dom \psi \mid \psi(x) \in \dom \phi \}$ and whose action on elements of its domain is, of course, the usual formula for composition.  Observe that for all $\phi, \psi \in \sF$, 
$T(\phi) T(\phi) = T(\phi \circ \psi)$.  Note that there may be elements of $[0,1]$ which are not in the domain of $\phi \circ \psi$ but which are mapped by the formula for $\phi \circ \psi$ into elements of $[0,1]$.  In particular, if $g,h \in G$ then $\beta_g \circ \beta_h$ is a restriction of $\beta_{gh}$.  (Simple examples where the restriction is proper can be obtained using $h = g^{-1}$; indeed, for some of these, 
$\beta_g \circ \beta_h$ is the ``empty'' transformation.)  It is now clear that $T(\beta_g) T(\beta_h) = T(\beta_g \circ \beta_h)$ is a restriction of 
$T(\beta_{gh})$.

If $Q$ is a $k$-fold product of isometries each of which is one of
$T_1, \dots, T_n$, then $QQ^*$ is a projection onto $L^2(J)$, where $J$ is some interval of length $1/n^k$.  Suitable sums of such projections yield projections associated with intervals of length 
$p/n^k, p \in \bbZ, 1 \leq p \leq n^k$.  If we left multiply a partial isometry $T(\beta_g)$ by such a projection, we obtain $T(\phi)$, where
$\phi \in \sF$.

Similarly, ``translations'' are obtained from products $QR^*$ where each of $Q$ and $R$ is a word of length $k$ in $T_1, \dots, T_n$.  More generally, 
for any $\phi \in \sF$, $T(\phi)$ can be written as a sum of words in
$T_1, \dots, T_n$ and their adjoints.

\begin{proof}[Proof of Theorem \ref{t:Onpca}]
For each $\phi \in \sF$, let 
$S(\phi) = \chi_{\ran \phi}U^g$, where $g$ is that element of $G$ for which $\phi$ is a restriction of $\beta_g$.  Note: to avoid extra notation, we view $\phi$ as acting on $[0,1]$ and $X$ simultaneously.  In the definition of
$S(\phi)$, $\ran \phi$ is a subset of $X$ -- a clopen subset of
$\ran \beta_g$.

For each $i=1,\dots,n$, let
\begin{displaymath}
  S_i = S\left(\beta_{\left(\frac {i-1}n,1 \right)} \right)
 = \chi_{[\frac {i-1}n, \frac in]} U^{(\frac {i-1}n,1)}.
\end{displaymath}
Observe that, for each $i$,
$\dom  \beta_{\left(\frac {i-1}n, 1\right)}$ is $[0,1]$ and
$\ran  \beta_{\left(\frac {i-1}n, 1\right)}$ is
$\left[\frac {i-1}n, \frac in \right]$, so that 
$\chi_{\left[\frac {i-1}n, \frac in\right]} \circ  
\beta_{\left(\frac {i-1}n, 1\right)} = \chi_{[0,1]}$.
Formula (\ref{inv}) above yields
\begin{displaymath}
  S_i^* = \chi_{\left[\frac {i-1}n, \frac in \right]} \circ
   \beta_{\left(\frac {i-1}n, 1\right)} 
  U^{\left(\frac {i-1}n, 1 \right)^{-1}} = \chi_{[0,1]}U^{(-(i-1),-1)}.
\end{displaymath}
We then have (using formula (\ref{mult}))
\begin{align*}
\sads{i} &= \chi_{[0,1]}U^{(-(i-1),-1)}
     \chi_{[\frac {i-1}n, \frac in]}
 U^{\left(\frac {i-1}n,1\right)} \\
&= \chi_{[0,1]} \cdot \chi_{[\frac {i-1}n, \frac in]} \circ
   \beta_{\left(\frac {i-1}n,1\right)} U^{(0,0)}\\
&= \chi_{[0,1]} U^{(0,0)} =I.
\end{align*}
Thus, each $S_i$ is an isometry.  Now
$\dom \beta_{(-(i-1),-1)}=[\frac {i-1}n, \frac in]$ and
$\ran \beta_{(-(i-1),-1)}= [0,1]$, so
$\chi_{[0,1]} \circ \beta_{(-(i-1),-1)}
 = \chi_{[\frac {i-1}n, \frac in]}$.
Therefore,

\begin{align*}
\ssad{i} 
&= \chi_{[\frac {i-1}n, \frac in]} U^{\left(\frac {i-1}n,1\right)}
\chi_{[0,1]}U^{(-(i-1),-1)} \\
&= \chi_{[\frac {i-1}n, \frac in]} \cdot \chi_{[0,1]} \circ
  \beta_{(-(i-1),-1)} U^{(0,0)} \\
&= \chi_{[\frac {i-1}n, \frac in]} U^{(0,0)}.
\end{align*}

This yields $\sum_i S_i S_i^* = \chi_{[0,1]} U^{(0,0)} =I$.
Thus, $S_1, \dots, S_n$ are $n$ isometries in
$C(X) \times_{\alpha} G$ with orthogonal range projections
whose sum is $I$.  This shows that
$C^*(S_1, \dots, S_n)$, the \cstar subalgebra of
$C(X) \times_{\alpha} G$ generated by $S_1, \dots, S_n$, is isomorphic
to the Cuntz algebra $O_n$.  To complete the proof we need merely show 
that $C^*(S_1, \dots, S_n)$ is all of $C(X) \times_{\alpha} G$.

With the help of formula (\ref{mult}) we can show that
$S(\phi \circ \psi) = S(\phi)S(\psi)$, for any $\phi,\psi \in \sF$.  This implies that if $T(\phi)$ is a sum of words in the $T^{\phantom{*}}_i$ and $T^*_i$, then $S(\phi)$ is the same sum of words in the
$S^{\phantom{*}}_i$ and $S^*_i$. From the known structure of $O_n$ described above, we deduce that, for any interval $J$,
$\chi_J U^{(0,0)} = \chi_J I$ is in the (non-closed) \star algebra generated by $S_1, \dots, S_n$ and that, for each $g \in G$, 
$\chi_{\ran \beta_g} U^g$ is also in this algebra.  The first of these two facts implies that $f U^{(0,0)} = fI \in C^*(S_1, \dots S_n)$, for all
$f \in C(X)$.  This, combined with the second fact, shows that
$f U^g \in C^*(S_1, \dots, S_n)$ for all 
$f \in C_0(\ran \beta_g)$.  It now follows immediately that
$C(X) \times_{\alpha} G \subseteq C^*(S_1, \dots S_n)$ and the proof is complete.
\end{proof}

\section{The connection between the partial action and the groupoid} 
\label{s:pagr}

As we have seen, the group $G$ and the partial action
studied in section \ref{s:On} are suggested by the Cuntz groupoid
model for $O_n$ and the standard representation of $O_n$ acting on
$L^2[0.,1]$. It is not surprising that in this situation the groupoid
and the partial action determine each other.  General theorems about
this connection can be found in \cite{MR2045419,apcp}.  The theorem in
\cite{apcp} does not apply directly to the partial action of $G$ on
$X$ and the Cuntz groupoid, since the partial action is not freely
actiong and the Cuntz groupoid is not principal.  In contrast to the
two references cited, the content of this section is naive: we simply
show directly how the partial action of $G$ determines the Cuntz
groupoid and vice versa. This connection does follow Theorem 5.1 
in \cite{apcp} in spirit.

When working with the Cuntz groupoid it is more convenient to regard
$X$ as the space of all sequences with entries from
$\{0,\dots,n-1\}$.  The correspondence
$(a_1, a_2, \dots) \leftrightarrow \sum \frac {a_i}{n^i}$ connects the
two different representations of the Cantor space $X$.  (Sequences
with a tail of $0$'s or with a tail of $n-1$'s correspond to $n$-adic
rationals with a superscript $-$ or $+$, of course.)  As a set, the
Cuntz groupoid $\sG$ is the set of all triples $(x,k,y)$ where
$x,y \in X$ and $x_{i+k} = y_i$ for all sufficiently large $i$.

\subsection{From partial action to groupoid}  This is the easy direction.
For each $(r,k) \in G$, define
\begin{displaymath}
\text{``graph''}(\beta_{(r,k)})  = \{ (x,k,y) \mid 
y \in \dom \beta_{(r,k)} \text{ and } x = \beta_{(r,k)}(y) \}.
\end{displaymath}
and let
\begin{displaymath}
  \sG = \bigcup_{(r,k) \in G} \text{``graph''}(\beta_{(r,k)}).
\end{displaymath}
The topology on $\sG$ is generated by all sets of the form
$\{(x,k,y) \mid x \in U \text{ and } y \in V\}$, where $V$ is a an
open subset of $ \dom \beta_{(r,k)}$ and
$U = \beta_{(r,k)}(V)$.
It is straightforward to check that $\sG$ is the Cuntz groupoid.

\subsection{From groupoid to partial action}  Here we need to define a 
$G$-valued cocycle $c$ so that the ``graph''($\beta_{(r,k)}$) turns
out to be $c^{-1}(r,k)$.
An element $(x,k,y) \in \sG$ clearly determines $k$; we need to
determine $r$ also so that 
$(x,k,y) \in \text{``graph''}\beta_{(r,k)}.$ 
For this purpose, we use
the standard
function which converts $k$-tuples from $\{0,1,\dots,n-1\}$ into 
$n$-adic rationals in $[0,1]$.  For 
$\lambda = (\lambda_1, \dots, \lambda_{|\lambda|})$, define
\begin{displaymath}
  s(\lambda) = \frac {\lambda_1}n + \dots +
  \frac {\lambda_{|\lambda|}}{n^{|\lambda|}}. 
\end{displaymath}
For two finite sequences $\lambda$ and $\mu$,
define
\begin{displaymath}
r(\lambda,\mu) = s(\lambda) - \frac {s(\mu)}{n^k}. 
\end{displaymath}

Finally define a cocycle $c \colon \sG \to G$ as follows:  
with $x=\lambda z$
and $y = \mu z$ as above, and $(x,k,y) \in \sG$, define
\begin{displaymath}
  c(x,k,y) = (r(\lambda,\mu), k).
\end{displaymath}
Note that on the set
$Z(\lambda,\mu) = \{(\lambda z, |\lambda|-|\mu|, \mu z \mid z \in X \}$
 $c$ has the constant value
$(r(\lambda,\mu),|\lambda|-|\mu|)$. 
Since these sets form a basis for the topology on $\sG$,
  $c$ is locally 
constant.

We do need to show that $c$ is a cocycle; for that purpose the
following concatenation forumla is useful.  For
finite strings $\eta$, $\zeta$,
\begin{displaymath}
  s(\eta \zeta) = s(\eta) + \frac {s(\zeta)}{n^{|\eta|}}.
\end{displaymath}

Let $(x,k,y) \in Z(\lambda,\mu)$ and 
$(y,l,z) \in Z(\gamma, \delta)$.  Assume that this is a composable
pair, which requires that one of $\mu$ and $\gamma$ extend the other.  Say, 
for the sake of argument, that $\gamma = \mu \gamma'$.  (The other
case can be handled in a similar way.)
We have
\begin{align*}
c(x,k,y) &= \left( s(\lambda) - \frac {s(\mu)}{n^k}, k \right) \text{ and} \\
c(y,l,z) &= \left( s(\gamma) - \frac {s(\delta)}{n^l}, l \right).
\end{align*}
The product of these two elements of $G$ is
\begin{displaymath}
  \left( s(\lambda) - \frac {s(\mu)}{n^k} 
+ \frac {s(\gamma)}{n^k} - \frac {s(\delta)}{n^{k+l}}, k+l \right).
\end{displaymath}
Now
\begin{align*}
c(x, k+l, z) &= \left(s(\lambda \gamma') -
\frac {s(\delta)}{n^{k+l}}, k+l \right) \\
&= \left( s(\lambda) + \frac {s(\gamma')}{n^{|\lambda|}} 
- \frac {s(\delta)}{n^{k+l}}, k+l \right).
\end{align*}
Therefore, we need to show that
\begin{displaymath}
  \frac 1{n^k} (s(\gamma) - s(\mu)) = 
\frac 1{n^{|\lambda|}} s(\gamma').
\end{displaymath}
Now, by the concatenation property,
\begin{displaymath}
  s(\gamma) = s(\mu) + \frac 1{n^{|\mu|}} s(\gamma');
\end{displaymath}
so
\begin{displaymath}
  s(\gamma) - s(\mu) = \frac 1{n^{|\mu|}}s(\gamma')
\end{displaymath}
and
\begin{displaymath}
  \frac 1{n^k}(s(\gamma) - s(\mu)) =
\frac 1{n^{k+|\mu|}} s(\gamma') =
\frac 1{n^{|\lambda|}} s(\gamma').
\end{displaymath}
This establishes the cocycle property for $c$.

Finally, we show that for each $(r,k) \in G$, $c^{-1}(r,k)$ is the
``graph'' of a partial homeomorphism on $X$. 
Let
\begin{align*}
&x = \lambda z, \\
&y = \mu z, \\
&k = |\lambda| - |\mu|, \\
&r = r(\lambda,\mu).
\end{align*}
We then have
\begin{align*}
x &\longleftrightarrow s(\lambda) + \frac 1{n^{|\lambda|}}
\sum \frac {z_i}{n^i}, \\
y &\longleftrightarrow s(\mu) + \frac 1{n^{|\mu|}}
\sum \frac {z_i}{n^i}
\end{align*}
and 
\begin{displaymath}
  r = s(\lambda) - \frac {s(\mu)}{n^k}.
\end{displaymath}
Calculate:
\begin{align*}
\beta_{(r,k)}(y) &= r + \frac y{n^k} \\
&= s(\lambda) - \frac {s(\mu)}{n^k} +
\frac 1{n^k}\left(s(\mu) + \frac 1{n^{|\mu|}}
\sum \frac {z_i}{n^i}   \right) \\
&= s(\lambda) + \frac 1{n^{k+|\mu|}} \sum \frac {z_i}{n^i} \\
&= s(\lambda) + \frac 1{n^{|\lambda|}} \sum \frac {z_i}{n^i} \\
&=x.
\end{align*}
This ties together the cocycle and the partial action.

\begin{remark}
Different pairs $\lambda$ and $\mu$ may yield the same value
for $r = r(\lambda, \mu)$.  But the conclusion
$\beta_{(r,k)}(y)=x$ is independent of $\lambda$ and $\mu$; it is
valid for any $(x,k,y)$ with $c(x,k,y) = (r,k)$.
\end{remark}

\begin{flushleft}
\begin{center}
\begin{flushright}

\end{flushright}
\end{center}
\end{flushleft}
\section{Subalgebras} \label{subalg}

One reason for focusing on the presentation of $O_n$ as a partial
action by $\bbQ_n \times_{\delta} \bbZ$ on $X$  is the possibility
that this may provide further insight into (non-self-adjoint)
subalgeras of $O_n$.  Some of these algebras have been studied in
\cite{hop_peters} and, in the more general context of graph \cstar
algebras, in \cite{hpp}.  In this section we describe a couple of
these subalgebras in terms of the partial crossed product
presentation.

The simplest way to obtain subalgebras is to generate them by
monomials associated with a subset of $G$.
If $P$ is a subset of $G$, let $B(P)$ be the closure in
$C(X) \times_{\alpha} G$ of the set of polynomials of the form
$\sum_{g \in P}f_g U^g$.  If $P$ is closed under multiplication, then
$B(P)$ is a subalgebra of $O_n$; if $P$ is closed under inversion then
$B(P)$ is self-adjoint.

So, for example,
if $P = \{(r,0) \mid r \in \bbQ_n \}$, then $B(P)$ is the
canonical $n^{\infty}$-UHF subalgebra of $O_n$. 
 If
$P = \{(r,0) \mid r \geq 0 \}$, then $B(P)$ is the refinement TAF
subalgebra of the canonical UHF subalgebra.  These subalgebras are not
very interesting as subalgebrsas of $O_n$.  For a more interesting
example, let 
$P = \{(r,k) \mid k>0 \text{ or } k = 0 \text{ and } r \geq 0 \}$.
Then $B(P)$ is a strongly maximal triangular subalgebra of $O_n$ whose
diagonal is the canonical masa in the canonical UHF subalgeba of
$O_n$.  Clearly, $B(P)$ is a superalgebra of the refinement TAF
algebra just mentioned.  Also, it is shown in \cite{hop_peters} that
$B(P)$ is semisimple.

The Volterra nest subalgebra $\sV$
of $O_n$ is another interesting example of a
non-self-adjoint subalgebra of $O_n$.  This was first introduced by
Power in \cite{MR86d:47057}, studied using groupoid techniques in
\cite{hop_peters} and extended to the graph \cstar algebra context in 
\cite{hpp}.  In the standard representation of $O_n$ acting on
$L^2[0,1]$, $\sV$ is just the intersection of $O_n$ and the usual
Volterra nest subalgebra acting on $L^2[0,1]$.  We now proceed to
identify $\sV$ in partial crossed product language.  The
reader is referred to \cite{hop_peters} for the groupoid theoretic
 description of the spectrum of $\sV$.

For each $n$-adic $r$ (greater than 0) in $X$, let
$p_r$ denote the characteristic function of $[0,r^-]$.  The elements 
$p_r U^0$ form a nest of projections in
 $C(X) \times_{\alpha} G = O_n$.  
This nest (with $0$ adjoined) is the Volterra nest in $O_n$.

Let $g$ be an element in $G$ and let $I$ be a clopen interval
contained in $\dom \beta_g$.  Since
$\chi_I \circ \beta_{g^{-1}} = \chi_{\beta_g(I)}$, we have
\[
\chi_{\ran \beta_g} U^g \chi_I U^0 = \chi_{\beta_g(I)} U^g.
\]
Suppose, further, that $\beta_g(t) \leq t$, for all $t \in I$.   Then 
$\chi_{\beta_g(I)} U^g$ leaves invariant all the Volterra projections
$p_r U^0$.  Indeed,
\[
(p_r U^0)^{\perp} \chi_{\beta_g(I)} U^g p_r U^0 =
p_r^{\perp} \chi_{\beta_g(I)} \cdot (p_r \circ \beta_{g^{-1}}) U^g,
\]
so we just have to show that
\[
p_r^{\perp}(t) \chi_{\beta_g(I)}(t) p_r(\beta_{g^{-1}}(t)) = 0,
\text{ for all } t.
\]

If $t \notin \beta_g(I)$, this expression is certainly 0. 
If  $\beta_{g^{-1}}(t) > r$, it is also zero.
So assume $t \in \beta_g(I)$ and $\beta_{g^{-1}}(t) \leq r$.
 Then $\beta_{g^{-1}}(t) \in I$  and
$t = \beta_g(\beta_{g^{-1}}(t)) \leq \beta_{g^{-1}}(t) \leq r$.
But then $p_r^{\perp}(t)=0$.  Thus the expression is zero for all
values of $t$.

For each $(r,k) \in G$, let $\Irk$ be the maximal clopen interval
contained in $\dom \brk$ with the property that
$\brk (t) \leq t$ for all $t \in \Irk$.  It is possible that
$\Irk$ is empty when $\dom \brk$ isn't; it is also possible that
$\Irk = \dom \brk$.  When $\Irk$ is a proper subinterval of
$\dom \brk$, there is an $n$-adic number $x$ such that
$\brk(x) = x$ (more accurately, $\brk(x^-) = x^-$ and
$\brk(x^+) = x^+$).  In this case, either $x^-$ will be the right hand
endpoint of $\Irk$ or $x^+$ will be the left hand endpoint.

\begin{remark}
It is not hard to work out when $\Irk = \dom \brk$.  Since
$\dom \brk = \emptyset$ when $r \geq 1$ we assume 
$r < 1$.  Then $\Irk = \dom \brk$ if, and only if, one of the
following two conditions holds:
\begin{enumerate}
\item $k \geq 0$ and $r \leq 0$,
\item $k < 0$ and $\ds r \leq 1 - \frac 1{n^k}$.
\end{enumerate}
\end{remark}

For each $(r,k) \in G$, define $\Trk = \chi_{\brk(\Irk)} \Urk$. 
 By the
comments above, $\Trk \in \sV$ for all $(r,k) \in G$.  In fact,
$\{\Trk \mid (r,k) \in G \}$ generates $\sV$.  This can be proven by
showing that
\begin{displaymath}
  \sigma(\sV) = \bigcup_{(r,k) \in G} \sigma(\Trk),
\end{displaymath}
where $\sigma$ denotes the spectrum in the groupoid.  The description
of $\sigma(\sV)$, which is given in \cite{hop_peters}, is a bit
involved and the verification of the equality above is
straightforward, so we omit the details.

\section{$M_k(O_n)$ as a partial crossed product} \label{s:MkOn}

$M_k(O_n)$ can also be written as a partial crossed product by a
partial action of an amenable group on an abelian algebra.  As
mentioned in the introduction, algebras more general than $M_k(O_n)$
are known to be partial crossed products (with non-amenable groups).
But the cosntruction described in this section may prove useful
in the study of subalgebras of $M_k(O_n)$.

Let $X$ be the $n$-adic Cantor space and 
$G = \bbQ_n \times_{\delta} \bbZ$, as in section \ref{s:On}.
Let $S_k = \{0, 1, \dots, k-1\}$.  The spectrum of the abelian algebra
used in the partial crossed product construction will be the Cartesian
product $Y = X \times S_k$ and the group will be the product group
$H = G \times \bbZ$.  The partial action is given by
\begin{displaymath}
  \beta_{(r,j,p)}(x,t) = \left(\frac x{n^j}+r,t+p \right).
\end{displaymath}
The domain of $\beta_{(r,j,p)}$ is
$\{(x,t) \in Y \mid \frac x{n^j}+r \in X \text{ and } t+p \in S_k\}$.
Let $\alpha$ be the partial action on $C(Y)$ induced by $\beta$ in the
usual way.

The partial action
$\beta$ is obviously built from the partial action of $G$ on $X$ in
section \ref{s:On}, which we shall now denote by $\beta^1$ and a
partial action $\beta^2$ of $\bbZ$ on $S_k$ given by
$\beta^2_p(t)=t+p$ on the obvious domain.  It is well known that 
$C(S_k) \times_{\alpha^2} \bbZ \cong M_n$, where $\alpha^2$ is
$\beta^2$ transferred to $C(S_k) \cong D_n$.

The fact that $C(Y) \times_{\alpha} H \cong M_k(O_n)$  follows
immediately from Theorem \ref{T:tensor} below.

In the special case when each $C_i$ is abelian, and so of the form
$C(X_i)$, and the action of each
 $G_i$ is topologically free, a shorter proof based on Theorem 2.6 of \cite{MR1905819}
is available.  These conditions are satisfied by the algebras and partial actions used for $M_k(O_n)$ and the shorter proof is sketched in a remark after the proof of Theorem \ref{T:tensor}.

\begin{theorem} \label{T:tensor}
For $i = 1,2$, let $C_i$ be a nuclear \cstar algebra, 
$G_i$ an amenable group, and $\alpha^i$ a partial action of $G_i$ on $C_i$.  Let $\alpha = \alpha^1 \otimes \alpha^2$ be the partial action of $G_1 \times G_2$ on $C_1 \otimes C_2$ defined on elementary tensors by
\begin{displaymath}
\alpha_{(g_1,g_2)}(c_1 \otimes c_2) 
= \alpha_{g_1}^1 \otimes \alpha_{g_2}^2 (c_1 \otimes c_2) 
= \alpha_{g_1}^1(c_1) \otimes \alpha_{g_2}^2(c_2), 
\quad c_1 \in \dom \alpha_{g_1}^1,
\;c_2 \in \dom \alpha_{g_2}^2.
\end{displaymath}
Let $A_1 = C_1 \times_{\alpha^1} G_1$,
$A_2 = C_2 \times_{\alpha^2} G_2$ and
$A_3 = (A_1 \otimes A_2) \times_{\alpha} (G_1 \times G_2)$.  Then
$A_1 \otimes A_2 \cong A_3$.
\end{theorem}

\begin{proof}
All the groups are amenable; by \cite[Prop. 4.2]{MR96i:46083} the
partial crossed produts are isomorphic to the reduced partial crossed
products.  In \cite{MR96i:46083} McClanahan shows how to construct
from a faithful representation of a \cstar algebra a faithful
representation of the reduced partial crossed product.  We use this
construction to prove the theorem.  First, we summarize (with the
suppression of details) McClanahan's construction.

Let $K$ be a group with a partial action $\alpha$ on a \cstar algebra
$C$.  Let $\pi \colon C \to \sB(\sH)$ be a representation.  For each
$g \in K$, there is a representation 
$\pi_g \colon \ran \alpha_g \to \sB(\sH)$ defined by
$\pi_g(a) = \pi(\alpha_{g^{-1}}(a))$.  This is extended to a
representation of all of $A$ (in a unique way) by means of an
approximate unit $u_{\lambda}$ in $\ran \alpha_g$:
$\pi'_{g}(a) = \lim \pi_g(u_{\lambda}a)$.  (The limit is taken in the
strong operator topology.)  A ``regular'' representation $\tilde{\pi}$
of $A$ acting on $\ell^2(K,\sH)$ is defined by
$(\tilde{\pi}(a) \xi)(g) = \pi'_g (a)\xi(g)$.
Let $\lambda$ denote the left regular representation of $K$ acting on
$\ell^2(K, \sH)$. Finally, define $\tilde{\pi} \times \lambda$ by
specifying the action on monomials:  
$\tilde{\pi} \times \lambda (f_g U^g) = \tilde{\pi}(f_g) \lambda_g$.

McClanahan shows that when $\pi$ is faithful, then 
$\tilde{\pi} \times \lambda$ is a faithful representation of the
reduced partial crossed product.  If the group is amenable, then
$\tilde{\pi} \times \lambda$ is a faithful representation of the
partial crossed product.

We shall now apply this to faithful representations $\pi_1$ and
$\pi_2$ of $C_1$ and $C_2$ acting on $\sH_1$ and $\sH_2$.  Let
$\pi = \pi_1 \otimes \pi_2$.  This is a faithful representation of
$C_1 \otimes C_2$ acting on
$\sB(\sH_1) \otimes \sB(\sH_2) \cong \sB( \sH_1 \otimes \sH_2)$, given
on monomials by 
$\pi (c_1 \otimes c_2) = \pi_1(c_1) \otimes \pi_2(c_2)$.

We claim that $\tilde{\pi} = \tilde{\pi}_1 \otimes \tilde{\pi}_2$.
(Both of these representations act on the Hilbert space
$\ell^2(G_1 \times G_2, \sH_1 \otimes \sH_2)$).  To verify this, begin
with $(g_1,g_2) \in G_1 \times G_2$.  Note that
$\ran \alpha_{(g_1,g_2)} = \ran \alpha^1_{g_1} \otimes \ran
  \alpha^2_{g_2}$.  
Apply $\pi_{(g_1,g_2)}$ to an elementary tensor $c_1 \otimes c_2$ in 
its domain:
\begin{align*}
\pi_{(g_1,g_2)}  (c_1 \otimes c_2) &=
\pi(\alpha_{(g_1,g_2)^{-1}}(c_1 \otimes c_2)) \\
&= \pi(\alpha_{g_1^{-1}} (c_1) \otimes \alpha_{g_2^{-1}}(c_2) ) \\
&= \pi_1(\alpha_{g_1^{-1}} (c_1)) \otimes \pi_2
(\alpha_{g_2^{-1}}(c_2)) \\
&= (\pi_1)_{g_1}(c_1) \otimes (\pi_2)_{g_2}(c_2) \\
&= \left((\pi_1)_{g_1} \otimes (\pi_2)_{g_2} \right) (c_1 \otimes c_2).
\end{align*}
It follows that
$\pi_{(g_1,g_2)} = (\pi_1)_{g_1} \otimes (\pi_2)_{g_2}$ on 
 $\ran \alpha_{(g_1,g_2)}$.

Take an approximate unit for $\ran \alpha_{(g_1,g_2)}$ consisting of
the tensor product of approximate units for $\ran \alpha^1_{g_1}$ and
$\ran \alpha^2_{g_2}$ and strong limits to get
\begin{displaymath}
  \pi'_{(g_1,g_2)} = (\pi'_1)_{g_1} \otimes (\pi'_2)_{g_2}.
\end{displaymath}
Now, for any elementary tensors $\xi_1 \otimes \xi_2$ in
$\ell^2(G_1, \sH_1) \otimes \ell^2(G_2, \sH_2)$ and 
$c_1 \otimes c_2$ in $C_1 \otimes C_2$, we have
\begin{align*}
(\tilde{\pi} (c_1 \otimes c_2) \xi_1 \otimes \xi_2) (g_1,g_2) &=
\pi'_{(g_1,g_2)} (c_1 \otimes c_2) (\xi_1 \otimes \xi_2) (g_1,g_2)
\\
&= (\pi'_1)_{g_1}(c_1) \xi_1(g_1) \otimes 
(\pi'_2)_{g_2}(c_2)) \xi_2(g_2) \\
&= (\tilde{\pi}_1 (c_1) \xi_1) g_1 \otimes 
(\tilde{\pi}_2 (c_2) \xi_2) g_2 \\
&= \left((\tilde{\pi}_1 \otimes \tilde{\pi}_2)
(c_1 \otimes c_2) (\xi_1 \otimes \xi_2)\right) (g_1, g_2).
\end{align*}
It follows that $\tilde{\pi}$ and 
$\tilde{\pi}_1 \otimes \tilde{\pi}_2$ agree everywhere.

Let $\lambda_1$ and $\lambda_2$ be the left regular representations of
$G_1$ and $G_2$ on $\ell^2(G_1, \sH_1)$ and $\ell^2(G_2, \sH_2)$.
Then $\lambda = \lambda_1 \otimes \lambda_2$ is the left regular
representation of $G_1 \times G_2$ acting on
$\ell^2 (G_1 \times G_2, \sH_1 \otimes \sH_2)$.
($\lambda(g_1,g_2) = \lambda_1(g_1) \otimes \lambda_2(g_2)$.)  We then
have
\begin{displaymath}
  \tilde{\pi} \times \lambda = (\tilde{\pi}_1 \times \lambda_1)
 \otimes (\tilde{\pi}_2 \times \lambda_2).
\end{displaymath}
Since $\tilde{\pi}_i \times \lambda_i$ is a faithful representation of
$C(X_i) \times_{\alpha^i} G_i$ for $i = 1,2$ and
$\tilde{\pi} \times \lambda$ is a faithful representation of
$C(X_1 \times X_2) \times_{\alpha^1 \otimes \alpha^2} (G_1 \times
G_2)$, the theorem follows.
\end{proof}

\begin{remark}
When each $C_i$ is abelian and each partial action is topologically free, Theorem 2.6 of \cite{MR1905819} permits a simpler proof of Theorem \ref{T:tensor}.  (``Topologically free'' is defined in 
terms of the dual partial action on the spectrum of the \cstar algebra: if $\beta$ is such a dual action on $X$, $\beta$ is 
\emph{topologically free} if, for each $g$ not the group 
identity, the set of fixed points of $\beta_g$ has empty interior).

Here is a sketch of the simpler proof.  For each $i$, let
$(\pi_i, v_i)$ be a covariant representation of 
$(C_i, G_i, \alpha_i)$ such that the corresponding representation
$\lambda_i$ of $A_i$ is faithful.  This implies that each $\pi_i$ is faithful.  Since the partial crossed product algebras are nuclear, $\lambda_1 \otimes \lambda_2$ is a faithful representation of 
$A_1 \otimes A_2$.  It is straightforward, albeit a little tedious, to check that $(\pi_1 \otimes \pi_2, v_1 \otimes v_2)$ is a covariant representation of 
$(C_1 \otimes C_2, G_1 \times G_2, \alpha_1 \otimes \alpha_2)$.  Let 
$\lambda_3$ be the corresponding representation of $A_3$.  Furthermore, 
if $a_iU^{g_i}$ are monomials in $A_i$, then another routine calculation shows that
$(\lambda_1 \otimes \lambda_2) (a_1 U^{g_1} \otimes a_2 U^{g_2}) =
\lambda_3 ( a_1 \otimes a_2 U^{(g_1,g_2)})$.  It follows that
$\lambda_1 \otimes \lambda_2 (A_1 \otimes A_2) = \lambda_3 (B_3)$.  Finally, $\pi_1 \otimes \pi_2$ is a faithful representation of
$C_1 \otimes C_2$, since each of $\pi_1$ and $\pi_2$ are faithful.  By Theorem 2.6 in \cite{MR1905819}, $\lambda_3$ is faithful.  It follows immediately that $A_1 \otimes A_2 \cong A_3$.  The proof in Theorem
\ref{T:tensor} above is the same argument applied to the regular representations.
\end{remark}

It is also possible to prove that 
$C(Y) \times_{\alpha} H \cong M_k(O_n)$ directly.   We shall sketch this briefly, omitting
most of the details.

Let $S_1, \dots, S_n$ be a set of generating isometries for $O_n$.
Let $E_{ij}$ be matrix units for $M_k$.  Then the following elements
of $M_k(O_n) \cong O_n \otimes M_k$ generate $M_k(O_n)$;
\begin{align*}
T_i &= S_i \otimes E_{11}, \quad i = 1, \dots, n-1, \\
T_n &= S_n \otimes E_{1n}, \\
R_i &= I \otimes E_{i+1,i}, \quad i = 1, \dots, k-1.
\end{align*}

$M_k(O_n)$ is, of course, a grpah \cstar algebra.  One simple graph
which yields $M_k(O_n)$ has $k$ vertices and $n$ loops.  At one
vertex, say $v_1$, there are $n-1$ edges which are self-loops (i.e.,
$v_1$ is both the source and the range for each of these edges.)  In
addition, there is one more cycle made up of $k$ edges which run
successively from $v_1$ to $v_2$ to $v_3$, etc. and finally from
$v_k$ back to $v_1$.  The generators $\{T_i, R_j\}$ satisfy the
Cuntz-Krieger relations for this graph.  ($T_1, \dots, T_{n-1}$
correspond to the $n-1$ self-loops; $R_1, \dots, R_{k-}, T_n$
correspond to the cycle of length $k$).

Any set of partial isometries which satisfies this set of
Cuntz-Krieger relations will generate an algebra isomorphic to
$M_k(O_n)$, since this algebra satisfies the Cuntz-Krieger uniqueness
theorem.  So, to identify $C(Y) \times_{\alpha} H$ as $M_k(O_n)$, we
need merely find in $C(Y) \times_{\alpha} H$ a set of partial
isometries which satisfy these Cuntz-Krieger relations and which also
generate $C(Y) \times_{\alpha} H$ as a \cstar algebra.

Here are a set of such generators:
\begin{align*}
T_i &= \chi_{\left[\frac {i-1}{n},\frac in\right] \times \{0\}} 
U^{\left(\frac {i-1}{n},1,0 \right)}, \quad i = 1, \dots, n-1, \\
T_n &= \chi_{\left[\frac {n-1}{n}, 1 \right] \times \{0\}}
U^{\left(\frac {n-1}{n}, 1, 1-k\right)}, \\
R_i &= \chi_{[0,1] \times \{i\}} U^{(0,0,1)}, \quad i = 1, \dots, k-1. 
\end{align*}

The projections which correspond to the vertices in the graph
described above are
\begin{displaymath}
  P_{v_i} = \chi_{[0,1] \times \{i-1\}} U^{(0,0,0)}, \quad i = 1,
  \dots, k.
\end{displaymath}
Routine calculations show that: 
\begin{enumerate}
\item For each $i = 1, \dots, n-1$, the initial space of $T_i$ is
$P_{v_1}$ and the final space is
$\chi_{\left[\frac {i-1}{n}, \frac in \right] \times \{0\}}
U^{(0,0,0)}$.
\item For each $i = 1, \dots, k-1$, the initial space of $R_i$ is
$P_{v_i}$ and the final space is $P_{v_{i+1}}$.
\item The initial space of $T_n$ is $P_{v_n}$ and the final space is
$\chi_{\left[\frac {n-1}{n}, \frac in \right] \times \{0\}}
U^{(0,0,0)}$.
\end{enumerate}
The $M_k(O_n)$--Cuntz-Krieger relations for these partial isometries
follow immediately.

The calculations that $\{T_i, R_j\}$ generate
$C(Y) \times_{\alpha} H$ as a \cstar algebra are a pain; they are best avoided entirely.

\section{The groupoid--partial action connection for $M_k(O_n)$.}
\label{s:grpa2}

As in section \ref{s:pagr}, we can construct a groupoid model for
$M_k(O_n)$ from the partial action of $H$ on $Y$, and vice versa.
Since $M_k(O_n)$ is a graph \cstar algebra, the usual groupoid model
is the one based on path space with shift equivalence.  For the
purposes of this section, however, it is better to use a slightly
different (but isomorphic) model.

With $X$ the $n$-adic Cantor space, we take
\begin{displaymath}
\sG = \{ ((\mu z, j), p, (\nu z, h)) \mid \mu z, \nu z \in X, j,h \in
S_k, |\mu| - |\nu| = p \}.
\end{displaymath}
The groupoid operations and topology are as expected.

Begin with the partial crossed product.  Let
$(r,j,p) \in H = (\bbQ_n \times_{\delta} \bbZ) \times \bbZ$.  
If $x \in X$ and  $i \in S_k$ satisfy
$\frac x{n^j}+r \in X$ and $i+p \in S_k$ then
$(x,i) \in \dom \beta_{(r,j,p)}$.  Let
$(y,q) = (\frac x{n^j}+r,i+p) = \beta_{(r,j,p)}(x,i)$.  Define
$\text{``graph''}\beta_{(r,j,p)}$ to be all such triples
$((x,i),p,(y,q))$ and we obtain
\begin{displaymath}
  \sG = \bigcup_{(r,j,p) \in H}\text{``graph''}\beta_{(r,j,p)}.
\end{displaymath}
This indicates the passage from the partial action to the groupoid.

For the other direction we need an $H$ valued cocycle defined on
$\sG$.  Simply define
\begin{displaymath}
  c((x,i), j, (y,q)) = (r,j,p),
\end{displaymath}
where $r = y-x = \frac x{n^j} -x$ and $p=q-i$.
This cocycle is locally constant; furthermore, given $(r,j,p) \in H$,
$c^{-1}(r,j,p)$ determines $\beta$: if
$((x,i), j, (y,q)) \in c^{-1}(r,j,p)$ define
\makebox{$(x,i) \underset{\beta_{(r,j,p)}}{\longrightarrow} (y,q)$} 
to get the
partial action.

%\bibliographystyle{amsplain}
%\bibliography{biblio}

\providecommand{\bysame}{\leavevmode\hbox to3em{\hrulefill}\thinspace}
\providecommand{\MR}{\relax\ifhmode\unskip\space\fi MR }
\providecommand{\MRhref}[2]{%
  \href{http://www.ams.org/mathscinet-getitem?mr=#1}{#2}
}
\providecommand{\href}[2]{#2}

\end{document}